\documentclass[oneside,11pt]{article}

\usepackage{amsmath}
\usepackage{amssymb}
\usepackage{amsthm}
\usepackage[alphabetic]{amsrefs}
\usepackage{array}
\usepackage[all]{xy}
\usepackage{fancyhdr}
\usepackage{euscript}
\usepackage{graphics,graphicx}
\usepackage{cancel}
\usepackage{fancybox}
\usepackage{verbatim}
\usepackage{tikz}
\usepackage{longtable}

\usepackage[tableposition=below]{caption}
\usepackage[colorlinks=false,hyperindex]{hyperref}

\hypersetup{
	colorlinks,
	citecolor=violet,
	linkcolor=violet,
	urlcolor=blue}

\DeclareFontFamily{OT1}{rsfs}{}
\DeclareFontShape{OT1}{rsfs}{n}{it}{<-> rsfs10}{}
\DeclareMathAlphabet{\mathscr}{OT1}{rsfs}{n}{it}

\newtheoremstyle%
{custom}%
{}%                         Space above
{}%													Space below
{}%													Body font
{}%                         Indent amount
{}%                         Theorem head font
{.}%                        Punctuation after heading
{ }%                        Space after heading
{\thmname{}%                Additional specifications for theorem head
\thmnumber{}%
\thmnote{\bfseries #3}}%

\newtheoremstyle%
{Theorem}%
{}%
{}%
{\itshape}%
{}%
{}%
{.}%
{ }%
{\thmname{\bfseries #1}%
\thmnumber{\;\bfseries #2}%
\thmnote{\;(\bfseries #3)}}%

\theoremstyle{Theorem}
\newtheorem{theorem}{Theorem}[section]

\theoremstyle{definition}

\newtheorem{remark}[theorem]{Remark}

\newtheorem{conjecture}[theorem]{Conjecture}

\newcommand{\QQ}{\mathbf{Q}}

\newcommand{\FF}{\mathbf{F}}
\newcommand{\Fq}{\mathbf{F}_q}

\newcommand{\PP}{\mathbf{P}}

\newcommand{\Res}{\operatorname{Res}}

\DeclareMathOperator{\Span}{Span}

\newcommand{\avlink}[1]{\href{http://www.lmfdb.org/Variety/Abelian/Fq/#1}{\textsf{#1}}}

\begin{document}

\title{Counterexamples to a Conjecture of Ahmadi and Shparlinski}
\author{Taylor Dupuy, Kiran Kedlaya, David Roe, Christelle Vincent}

\maketitle

\begin{abstract}
Ahmadi-Shparlinski conjectured that every ordinary, geometrically simple Jacobian over a finite field has maximal angle rank. 
Using the L-Functions and Modular Forms Database, we provide two counterexamples to this conjecture in dimension 4.
\end{abstract}

%%%%%%%%%%%%%%%%%%%%%%%%%%%%%%%%%%%
\section{Introduction}
%%%%%%%%%%%%%%%%%%%%%%%%%%%%%%%%%%%

The following is a conjecture of Ahmadi--Shparlinski:
\begin{conjecture}[{\cite[\S 5]{Ahmadi2010}}] \label{conj:Ahmadi-Shparlinski}
	Every ordinary, geometrically simple Jacobian over a finite field has maximal angle rank. 
\end{conjecture}
In this paper we report that this conjecture is false. 
Our work used the L-Functions and Modular Forms Database (LMFDB), specifically its database of abelian varieties over finite fields which can be found here:
\begin{center}
\url{https://www.lmfdb.org/Variety/Abelian/Fq/}.
\end{center}
Documentation, further conjectures, and interesting statistics are reported in \cite{DKRV}. 

Apart from the counterexamples of Section \ref{sec:Ahmadi-Shparlinski}, this article briefly recalls the notion of angle rank in Section \ref{sec:anglerank}, presents how geometric simplicity is computed in the LMFDB (Section \ref{sec:simple}) and how Jacobians are tested for in the LMFDB (Section \ref{sec:Jacobians}). Following this, we describe our search (Section \ref{sec:search}), and because angle ranks are computed numerically in the LMFDB, we provide a proof of the computation of the angle rank for both examples in Section \ref{sec:proof}.
Readers can verify these counterexamples themselves using the code provided at
\begin{center}
\url{https://github.com/LMFDB/abvar-fq/},
\end{center}
which uses Sage \cite{sage}, PARI \cite{pari} and Magma \cite{magma}. 
Finally, the two counterexamples are in Section \ref{sec:Ahmadi-Shparlinski}.
We remark that in addition to providing counterexamples to the conjecture, we give two new methods for algebraically certifying angle ranks (as remarked above, as of January 2020, in the LMFDB the angle ranks are computed numerically using an LLL algorithm).

\begin{remark}
We began searching for counterexamples to Conjecture \ref{conj:Ahmadi-Shparlinski} since it is incompatible with the Shankar-Tsimerman conjecture \cite[Conj.~2.5]{Shankar2018}
which states that every simple abelian fourfold over $\overline{\FF}_p$ is isogenous to a Jacobian; since angle rank, ordinarity, and geometric simplicity are preserved under base change, this would imply that every simple abelian fourfold has maximal angle rank.
\end{remark}

\subsection*{Acknowledgements}
The authors began this project during the semester program ``Computational aspects of the Langlands program'' held at ICERM in fall 2015.
Subsequent workshops in support of the project were sponsored by the American Institute of Mathematics
and the Simons Foundation. 
In addition Dupuy was partially supported by the European Research Council under the European Unions Seventh Framework
Programme (FP7/2007-2013) / ERC Grant agreement no. 291111/ MODAG while working on this project.
Kedlaya was supported by NSF (grants DMS-1501214, DMS-1802161), IAS (visiting professorship 2018--2019), UCSD (Warschawski Professorship), and a Guggenheim Fellowship (fall 2015).  
Roe was supported by Simons Foundation grant 550033.
Vincent was partially supported by NSF (grant DMS-1802323).
%----------------------------------------------------
\section{Frobenius Angle Rank} \label{sec:anglerank}
%----------------------------------------------------

This section very briefly presents the definition of the angle rank of an abelian variety defined over a finite field; for more context and a longer discussion the reader is directed to \cite{DKRV}*{\S2.6, 3.8} and \cite{Dupuy2020b}.

For $A$ an abelian variety of dimension $g$ with $L$-polynomial $L(T) = \prod_{i=1}^{2g}(1-\alpha_i T)$, the \emph{angle rank} of $A$ is the quantity
\[
\delta(A) = \dim_{\QQ}( \Span_{\QQ}(\lbrace \arg(\alpha_i): 1\leq i \leq 2g  \rbrace \cup \lbrace \pi \rbrace  ) ) - 1 
\in \{0,\dots,g\}.
\]

The angle rank detects multiplicative relations among the roots of $L$; these are closely related to exceptional Hodge classes on powers of $A$.
For example, by a theorem of Zarhin \cite[Theorem~3.4.3]{Zarhin1994}, $\delta(A) = g$ if and only if there are no exceptional Hodge classes on any power of $A$. At the other extreme, $\delta(A) = 0$ if and only if $A$ is supersingular.

\section{Geometric Simplicity}\label{sec:simple}

To check that an ordinary isogeny class defined over $\FF_q$ is geometrically simple, we compute the tensor square of the $L$-polynomial,
note that it has no nontrivial cyclotomic factors, and 
apply \cite{CMSV2019}*{Lemma 7.2.7 (b)} to deduce that
 all of the geometric endomorphisms are defined over $\Fq$; in particular,
 each simple isogeny factor is geometrically simple.
In the counterexamples presented in this article, ordinarity and simplicity over $\Fq$ follow from the irreducibility of the $L$-polynomial.
For a discussion of the geometric endomorphism algebra in more generality, see \cite{DKRV}*{\S3.5}.

\section{Searching for Jacobians}\label{sec:Jacobians}

The current version of the LMFDB contains substantial data about whether isogeny classes contain Jacobians up to dimension 3 (see \cite[\S3.7]{DKRV}). However, for this paper we need data in dimension 4, for which the LMFDB currently contains only negative results
(e.g., a given isogeny class may fail to contain a Jacobian because it contains no principally polarizable
variety, or because it corresponds to an impossible sequence of point counts on a curve). We thus cannot use the LMFDB alone to certify that a given isogeny class contains a Jacobian.

We instead take the approach of constructing curves, computing their zeta functions, and matching their numerators to the Weil polynomials contained in the LMFDB.
To construct the curves, we note that in genus 4, every nonhyperelliptic curve is the intersection of a quadric and a cubic in $\PP^3$.
In practice, over $\FF_2$, this allows us to make an exhaustive search over both hyperelliptic and nonhyperelliptic curves,
using Magma to compute zeta functions.
(In the nonhyperelliptic case, we limit the options for the quadric as in \cite{Savitt2003}.)
Over $\FF_3$ and $\FF_5$, we enumerate only over hyperelliptic curves, using Sage to compute zeta functions.
This was enough to find the two counterexamples presented here.

\section{Results of the Search}\label{sec:search}
The Ahmadi-Shparlinski conjecture is a theorem in dimension 2, even without the ordinary condition \cite[Theorem~2]{Ahmadi2010}.
It is also a theorem in dimension 3, but this time it requires the ordinary condition
\cite[Theorem~1.1]{Zarhin2015}. 
We verified the consistency of these results with the LMFDB database.

In dimension 4 over $\FF_2$ there are 52 isogeny classes of ordinary, geometrically simple abelian varieties with angle rank at most 3 (in fact they are all equal to 3). 
Searching through curves, we found 620 distinct zeta functions, none of which occur among the previous list of 52 isogeny classes. 
Therefore there are no such 4-dimensional Jacobians over $\FF_2$, and the conjecture holds in this case. 

By contrast, the conjecture fails in dimension 4 over $\FF_3$ and $\FF_5$, as shown by the examples presented in Section \ref{sec:Ahmadi-Shparlinski}.

%---------------------------------------------------
\section{Algebraic Certification of Angle Ranks}\label{sec:proof}
%---------------------------------------------------
We now describe the procedure we used to compute a rigorous upper bound on $\delta(A)$. 
See \cite[\S 3.8]{DKRV} for an alternate approach that also gives a rigorous lower bound (which we do not need here).

 Let $P(T)$ be the Weil polynomial of an abelian variety $A$ over $\Fq$.
Fix a precision $\rho = \sigma^2$ (we default to $\rho = 625$). 

\begin{enumerate}
\item Compute the roots $\{\alpha_i\}$ of $P(T)$ that have positive imaginary part, in a complex field $C$ of precision $\rho$.
Set $t_i = \arg(\alpha_i) / \pi$, where $\arg$ is the principal branch of the logarithm.
Throw away duplicates, obtaining $0 < t_1 < \dots < t_m < 1$.
\item Use LLL to find independent integer relations $R_1, \dots, R_s$ among $\{t_1, \dots, t_m, 1\}$.
A relation is considered spurious if all coefficients are larger than $2^\sigma$,
and we interrupt the computation if some value is larger than $2^\sigma$ but others are not
(this did not happen for any isogeny class in the database).
The \emph{numerical angle rank} is $m - s$.
\item Find a number field $K$ in which $P(T)$ splits completely.
Choose an embedding $\iota \colon K \hookrightarrow C$ and let $\beta_1, \dots, \beta_m$
be the roots of $P(T)$ in $K$ with $0 < \arg(\iota(\beta_1)) < \dots < \arg(\iota(\beta_m)) < \pi$.
(In other words, the root $\beta_i$ has argument approximated by $t_i$ above.)
\item The roots $\beta_1, \dots, \beta_m$ together with the relations $R_1, \dots, R_s$ provide a certificate that the angle rank of $A$
is at most $m - s$.  One can check using exact arithmetic in $K$ that a relation $R_i = (c_1, \dots, c_{m+1})$ holds by confirming that
\[
(-q)^{c_{m+1}} \prod_{i=1}^m \beta_i^{c_i} = 1.
\]
\end{enumerate}

\begin{remark}
The upper bound $m-s$ can only fail to be sharp if at some point we discarded a relation as spurious when it was real.  Such a relation would have all coefficients larger than $2^\sigma$.
\end{remark}

\begin{remark}
Let $P(T) = \prod_{i=1}^{2g} (T-\alpha_i)$ be a Weil polynomial,
with roots ordered so that $\alpha_i\alpha_{i+g} = q$ (where the indices are taken modulo $2g$).
This remark explains an alternative method, using resultants and ``taking cyclotomic parts,'' to verify 
the existence of a relation of the form
\begin{equation}\label{eqn:nontrivial}
   \alpha_{i_1}^{e_1} \alpha_{i_2}^{e_2} \cdots \alpha_{i_j}^{e_j} = \zeta q^{j/2},
\end{equation}
where $e_1,\dots,e_j$ are specified positive integers and $\zeta$ is an unspecified root of unity.  
Such a relation is guaranteed to be nontrivial (i.e., not a consequence of the known relations
$\alpha_i \alpha_{i+g} = q$) provided that $i_1,\dots,i_j$ are pairwise distinct mod $g$;
when this occurs, the existence of the relation implies that the angle rank is not maximal. Note however
that as presented, this method cannot always certify a sharp upper bound on the angle rank.

Before presenting this method, we present three preliminary facts which we will need and take for granted:
\begin{enumerate}
	\item If $F(T)$ and $G(T)$ are any two polynomials, the polynomial 
	\begin{equation*}
	H(T)=\Res_S(F(S),G(T/S)S^{\deg G})
	\end{equation*}
	is a polynomial whose roots are products of the roots of $F$ and $G$.

	\item If $A$ is an abelian variety over $\FF_q$ of dimension $g$ with characteristic polynomial $P(T) = P_{A_{\FF_q}}(T) = \prod_{i=1}^{2g}(T-\alpha_i)$, then for any positive integer $n$,
	$P_{A_{\FF_{q^n}}}(T) = \prod_{i=1}^{2g}(T-\alpha_i^n)$.
	One can show that there is a formula for this polynomial in terms of resultants given by $P_{A_{\FF_{q^n}}}(T) = \Res_S(P(S), S^n-T).$
	\item If $F(T)$ is a polynomial with integer coefficients, then it factors as $F(T) = C(T)G(T)$ where $C(T)$ is a cyclotomic polynomial and $G(T)$ has no cyclotomic factors. The computation of the cyclotomic part $C(T)$ can be done efficiently using an algorithm described in \cite{BS02}; this is implemented in Sage by the function $\texttt{cyclotomic\_part()}$ called on $F(T)$, which returns $C(T)$.
\end{enumerate}

With these facts granted, note first that we are free to make the test after performing a base change as in the second point above; we may thus assume without loss of generality that $q$ is a perfect square, so that
$\overline{P}(T)=P(q^{1/2}T)$ is a root-unitary polynomial with rational coefficients. 
Set $\widetilde{\alpha}_i = q^{-i/2} \alpha_i$, so that the roots of $\overline{P}(T)$
are $\widetilde{\alpha}_1,\dots,\widetilde{\alpha}_{2g}$.

Using the first and second point, we produce the polynomial  
 $$ Q(T) = \prod_{i_1,i_2,\ldots,i_j}(T-\widetilde{\alpha}_{i_1}^{e_1} \widetilde{\alpha}_{i_2}^{e_2} \cdots \widetilde{\alpha}_{i_j}^{e_j}),$$
 where the product is taken over all tuples $(i_1,i_2,\ldots,i_j)$ in which $1 \leq i_k \leq 2g$. 
 We then use the third point to compute the cyclotomic part of $Q(T)$,
 and compare its degree with that of the ``trivial cyclotomic factor'' coming from the relations $\alpha_i\alpha_{i+g}=q$. 
(For example when $j=2$ and $e_1 = e_2 = 1$, $Q(T)$ is divisible by $(T-1)^g$.)
If these degrees disagree, this implies the existence of a nontrivial relation. (However, this relation might be ``partially trivial'' and thus involve fewer than $j$ distinct roots.)

In practice, we predict $j$ and $e_1,\dots,e_j$ in the relation \eqref{eqn:nontrivial} using LLL, and then verify the existence of a relation as we have just described, by showing that its cyclotomic part is greater than predicted by the trivial relations.\end{remark}

%---------------------------------------------------
\section{Counterexamples}	\label{sec:Ahmadi-Shparlinski}
%---------------------------------------------------

\subsection{A Counterexample when $g=4$ and $q=3$}
Let $C$ be the hyperelliptic curve over $\FF_3$ given by
\[
	y^2 = x^9 + x^8 + x^7 + 2x^5 + x.
\]
Then
\[
	L(C/\FF_3, T) = 
	1 - T + 2T^2 - 4T^3 - 2T^4 - 12T^5 + 18T^6 - 27T^7 + 81T^8,
\]
so the Jacobian $A$ of $C$ belongs to isogeny class \avlink{4.3.ab\_c\_ae\_ac}, which is ordinary and geometrically simple.
We compute the minimal splitting field of this polynomial, which has degree $48$ over $\QQ$, and fix an ordering of the roots of $P_A(T)$ such that $\beta_i\beta_{i+4}=3$, where the indices are taken modulo $8$.
Using the method of Section \ref{sec:proof}, we show that these roots satisfy the nontrivial relation
\[
\beta_1 \beta_3 \beta_4 = -3 \beta_2.
\]
The angle rank of $A$ is thus at most $3$, and is equal to $3$ unless there is a relation with exponents all larger than $2^{25}$.

\subsection{A Counterexample when $g=4$ and $q=5$}
Similarly, let $C$ be the hyperelliptic curve over $\FF_5$ given by 
\[
	y^2 = x^9 + x^6 + 2x^5 + x.
\]
Then
\[
	L(C/\FF_5,T) = 
	1-T+2T^2-4T^3+16T^4-20T^5+50T^6-125T^7+625T^8,
\]
so the Jacobian $A$ of $C$ belongs to isogeny class \avlink{4.5.ab\_c\_ae\_q}.
Again, $A$ is ordinary, geometrically simple, and has angle rank bounded above by $3$. 
If once again we order the roots of $P_A(T)$ such that $\beta_i\beta_{i+4}=3$, where the indices are taken modulo $8$, there is now a nontrivial relation of the form
\[
\beta_1 \beta_4 = \beta_2 \beta_3.
\]

\bibliographystyle{alpha}
\bibliography{weil}
\end{document}